\documentclass[]{amsart}
\usepackage{latexsym,amssymb,amsopn,amsbsy,amsthm,amsmath,amscd}

\newcommand{\g}{\mathfrak{g}}
\newcommand{\h}{\mathfrak{h}}
\newcommand{\p}{\mathfrak{p}}

\newcommand{\bbS}{\mathbb{S}}

\newcommand{\dirac}{\mbox{$\not\negthinspace\partial$}}

\newcommand{\Spinc}{\mathrm{Spin}^{c}}

\DeclareMathOperator{\ad}{ad}

\DeclareMathOperator{\Cl}{Cl}

\DeclareMathOperator{\ads}
{\widetilde{a\thickspace\thinspace}\negthickspace\negthinspace d}

\DeclareMathOperator{\End}{End} \DeclareMathOperator{\Hom}{Hom}
\DeclareMathOperator{\Ker}{Ker}

\DeclareMathOperator{\Index}{Index}

\theoremstyle{plain}

\newtheorem{theorem}{Theorem}
\theoremstyle{remark}

\begin{document}

\title{Harmonic Spinors on Homogeneous Spaces}

\date{February 24, 2000}

\author{Gregory~D. Landweber}

\address{Microsoft Research\\
         One Microsoft Way\\
         Redmond, WA  98052}

\email{gregland@microsoft.com}



\begin{abstract}
    Let $G$ be a compact, semi-simple Lie group and $H$ a
    maximal rank reductive subgroup. The irreducible
    representations of $G$ can be constructed as spaces of
    harmonic spinors with respect to a Dirac operator on
    the homogeneous space $G/H$ twisted by bundles associated to
    the irreducible, possibly projective, representations of $H$.
    Here, we give a quick proof of this result, computing the
    index and kernel of this twisted Dirac operator using a
    homogeneous version of the Weyl character formula noted by
    Gross, Kostant, Ramond, and Sternberg, as well as recent work
    of Kostant regarding an algebraic version of this Dirac operator.
\end{abstract}

\maketitle

\section{Introduction}

Recently Gross, Kostant, Ramond, and Sternberg discovered a
generalization of the Weyl character formula, describing the
character of any irreducible representation $V_{\lambda}$ of a
semi-simple Lie algebra $\g$ in terms of the characters of a
multiplet of irreducible representations $U_{\mu}$ of a maximal
rank reductive Lie subalgebra $\h$. Kostant then showed that these
$\h$-representations $U_{\mu}$ could be constructed as the kernel
of a certain Dirac operator acting on $V_{\lambda}\otimes\bbS$,
where $\bbS$ is the complex spin representation associated to the
complement of $\h$ in $\g$. These results are summarized in
Section \ref{section:background}, and a thorough exposition of the
subject can be found in \cite{St}.

After a brief review of homogeneous differential operators in
Section \ref{section:homogeneous}, we turn this argument on its
head in Section \ref{section:newstuff}. Let $G$ be the compact,
simply connected Lie group with Lie algebra $\g$, and let $H$ be
the Lie subgroup with Lie algebra $\h$. If the homogeneous space
$G/H$ is a spin manifold, then we show that the index of a
standard geometric Dirac operator on $G/H$ twisted by the
homogeneous vector bundle induced by a representation $U_{\mu}$ of
$H$ is, up to sign, the corresponding representation $V_{\lambda}$
of $G$. Furthermore, using the geometric version of Kostant's
Dirac operator, we explicitly construct this representation
$V_{\lambda}$ as the space of twisted harmonic spinors. If $G/H$
is not spin, or if the twisted spinors do not come from a suitable
$\Spinc$ structure on $G/H$, then we instead carry out these
constructions using the analogous $\h$-equivariant operators
upstairs on the Lie group $G$.

\section{The Homogeneous Weyl Formula}
\label{section:background}

Let $\g$ be a compact, semi-simple Lie algebra, and let $\h$ be a
reductive Lie subalgebra of maximal rank in $\g$. Since $\h$ has
the same rank as $\g$, any Cartan subalgebra of $\h$ is likewise a
Cartan subalgebra of $\g$, and the roots of $\h$ are then a subset
of the roots of $\g$. The Weyl group $W_{\g}$ of $\g$ acts simply
transitively on the Weyl chambers for $\g$, each of which is
contained inside a Weyl chamber for $\h$. Choosing a set of
positive roots for $\g$ also determines a positive root system for
$\h$, and we define $C\subset W_{\g}$ to be the subset of elements
that map the positive Weyl chamber for $\g$ into the positive Weyl
chamber for $\h$.

Let $\rho_{\g}$ and $\rho_{\h}$ denote half the sum of the
positive roots of $\g$ and $\h$ respectively. Equivalently, $\rho$
is also the sum of the fundamental weights, which lie on the
boundary of the positive Weyl chamber. For any dominant weight
$\lambda$ of $\g$, the weight $\lambda + \rho_{\g}$ lies in the
interior of the positive Weyl chamber for $\g$. For any $c\in C$,
the weight $c(\lambda+\rho_{\g})$ then lies in the interior of the
positive Weyl chamber for $\h$, and so
\begin{equation*}
    c\bullet \lambda := c\,( \lambda + \rho_{\g} ) - \rho_{\h}
\end{equation*}
is a dominant weight for $\h$. The $\rho$-shift insures that each
of the weights $c\bullet\lambda$ for $c\in C$ is distinct. Note
that although every dominant weight of $\g$ corresponds to a
distinct multiplet of dominant weights of $\h$, not every dominant
weight of $\h$ corresponds to a dominant weight of $\g$. In
particular, if $\mu$ is a dominant weight of $\h$ such that $\mu +
\rho_{\h}$ lies on the boundary of a Weyl chamber for $\g$, then
$\mu$ is not of the form $c\bullet\lambda$ for any dominant weight
$\lambda$ of $\g$. Such orphan weights will behave as exceptional
cases in the results of Section \ref{section:newstuff}.

Putting an $\ad$-invariant inner product on $\g$, let $\p$ denote
the orthogonal complement to $\h$ in $\g$. The adjoint action of
$\g$ then restricts to give an orthogonal action of $\h$ on $\p$,
and this action lifts to the complex spin representation $\bbS$
associated to $\p$. Furthermore, since $\h$ is of maximal rank in
$\g$, the complement $\p$ is even dimensional, and so the spin
representation decomposes as the sum $\bbS =
\bbS^{+}\oplus\bbS^{-}$ of two distinct half-spin representations.
To specify the sign convention, note that the weight space of
$\bbS$ with highest weight $\rho_{\g} - \rho_{\h}$ is one
dimensional, and take $\bbS^{+}$ to be the half-spin
representation containing that highest weight space.

In \cite{GKRS}, Gross, Kostant, Ramond, and Sternberg prove

\begin{theorem}[Homogeneous Weyl Formula]
    Let $V_{\lambda}$ be the irreducible representation of $\g$
    with highest weight $\lambda$. The following identity
    holds in the representation ring $R(\h)$:
    \begin{equation}\label{eq:homogeneous-weyl}
        V_{\lambda} \otimes \bbS^{+} -
        V_{\lambda} \otimes \bbS^{-}
        = \sum_{c\in C}(-1)^{c}\,U_{c\bullet\lambda},
    \end{equation}
    where $V_{\lambda}$ on the left side is viewed as a
    representation of $\h$ by restriction, and
    $U_{c\bullet\lambda}$ denotes the irreducible representation
    of $\h$ with highest weight $c\bullet\lambda$.
\end{theorem}

If $\h$ is a Cartan subalgebra $\mathfrak{t}$ of $\g$, then $C$ is
the full Weyl group $W_{\g}$, which acts on weights as
$w\bullet\lambda = w\,(\lambda+\rho_{\g})$, and
(\ref{eq:homogeneous-weyl}) becomes the familiar Weyl character
formula
\begin{equation}\label{eq:weyl}
    \chi(V_{\lambda}) \otimes
    \bigl( \, \bbS_{\g/\mathfrak{t}}^{+} -
           \bbS_{\g/\mathfrak{t}}^{-} \, \bigr)
    = \sum_{w\in W_{\g}} (-1)^{w}\, w (e^{\lambda + \rho_{\g}} ).
\end{equation}
The general form (\ref{eq:homogeneous-weyl}) of this identity can
be derived from the Weyl character formula by dividing both sides
of (\ref{eq:weyl}) by the character of the virtual spin
representation $\bbS_{\h/\mathfrak{t}}^{+} -
\bbS_{\h/\mathfrak{t}}^{-}$ associated to the complement of
$\mathfrak{t}$ in $\h$.

In \cite{K}, Kostant constructs a Dirac operator
$\dirac:V_{\lambda}\otimes\bbS^{+} \rightarrow
V_{\lambda}\otimes\bbS^{-}$. Since the domain and range are finite
dimensional, the $\h$-index of any such operator is automatically
given by (\ref{eq:homogeneous-weyl}). However, Kostant's Dirac
operator is unique in that it also respects the sign decomposition
given by the right side of (\ref{eq:homogeneous-weyl}), satisfying
\begin{equation}\label{eq:kostant-kernel}
    \Ker\dirac   = \sum_{(-1)^{c}=+1} U_{c\bullet\lambda},
    \quad
    \Ker\dirac^{\ast} = \sum_{(-1)^{c}=-1} U_{c\bullet\lambda}.
\end{equation}
This Dirac operator is formally self-adjoint, so the adjoint of
the Dirac operator $\dirac^{\ast}:V_{\lambda}\otimes\bbS^{-}
\rightarrow V_{\lambda}\otimes\bbS^{+}$ can be viewed as the same
operator acting on the opposite half-spin representation. This
operator thus provides a mechanism by which to extract the
multiplet of $\h$-representations $U_{c\bullet\lambda}$ directly
from the associated $\g$-representation $V_{\lambda}$.

In its most abstract form, Kostant's Dirac operator can be viewed
as an element of the non-abelian Weil algebra $U(\g) \otimes
\Cl(\p)$ (see \cite{AM}), where $U(\g)$ is the universal
enveloping algebra of $\g$, and $\Cl(\p)$ is the Clifford algebra
of $\p$. Choosing a basis $\{X_{i}\}$ of $\p$ and letting
$\{X_{i}^{\ast}\}$ be the dual basis satisfying $\langle
X_{i},X_{j}^{\ast}\rangle = \delta_{ij}$, Kostant defines his
Dirac operator to be the element
\begin{equation}\label{eq:kostantdirac}
    \dirac := \sum_{i} X_{i}\otimes X_{i}^{\ast}
            + 1 \otimes v,
\end{equation}
where $v\in\Cl(\p)$ is the image of the fundamental 3-form
$\omega\in\Lambda^{3}(\p^{\ast})$,
\begin{equation}\label{eq:threeform}
    \omega(X,Y,Z) = \langle X,[Y,Z] \rangle,
\end{equation}
under the Chevalley identification
$\Lambda^{\ast}(\p^{\ast})\rightarrow\Cl(\p)$. Now, any
representation $r$ of $\g$ on $V_{\lambda}$ extends to a
homomorphism $r:U(\g)\rightarrow\End(V_{\lambda})$, and the
Clifford action on the spin representation yields a homomorphism
$c:\Cl(\p)\rightarrow\End(\bbS)$ with the odd part of the Clifford
algebra interchanging $\bbS^{+}$ and $\bbS^{-}$. Combining these
maps gives a representation of the non-abelian Weil algebra on the
tensor product,
\begin{equation*}
    r \otimes c : U(\g)\otimes\Cl(\p) \rightarrow
                  \End(V_{\lambda}\otimes\bbS),
\end{equation*}
which takes the Dirac element $\dirac$ to an operator in
$\Hom(V_{\lambda}\otimes\bbS^{+},V_{\lambda}\otimes\bbS^{-})$ as
desired.

To compute the kernel of this Dirac operator, Kostant expressed
its square as a sum of quadratic Casimir operators for the Lie
algebras $\g$ and $\h$. When restricted to a subspace of
$V_{\lambda}\otimes\bbS$ transforming like the representation
$U_{\mu}$ under the diagonal right $\h$-action, the square of the
Dirac operator becomes
\begin{equation}\label{eq:dsquared}
    \dirac^{2}|_{U_{\mu}} = \| \lambda + \rho_{\g} \|^{2} -
                            \| \mu + \rho_{\h} \|^{2}.
\end{equation}
Since the Weyl group acts by isometries, all of the weights $\mu$
of $\h$ satisfying $\mu + \rho_{\h} = c(\lambda + \rho_{\g})$ must
have the same $\rho$-shifted norm as $\lambda$.  In fact, the
weights $\mu = c\bullet\lambda$ are precisely those for which the
expression (\ref{eq:dsquared}) vanishes. Furthermore, Kostant
showed that each of the representations $U_{c\bullet\lambda}$
occurs exactly once in the decomposition of
$V_{\lambda}\otimes\bbS$, and it follows that the kernel of the
Dirac operator must be of the form (\ref{eq:kostant-kernel}).

\section{Homogeneous Differential Operators}
\label{section:homogeneous}

Looking at the discussion of the previous section from a
geometric, rather than algebraic viewpoint, let $G$ be a compact,
semi-simple Lie group, and let $H$ be a reductive Lie subgroup of
maximal rank. Without loss of generality, also assume that $G$ is
simply connected. The Hilbert space $L^{2}(G)$ of functions on $G$
admits a $G\times G$ action, with left and right components
\begin{equation}\label{eq:ggaction}
    l(h) f : g \mapsto f(h^{-1}g), \quad
    r(h) f : g \mapsto f(gh)
\end{equation}
for any function $f\in L^{2}(G)$ and elements $g,h\in G$. At the
Lie algebra level, the infinitesimal left and right actions on
functions are given by differentiation with respect to the
right-invariant and left-invariant vector fields respectively.
More precisely, for any $X\in\g$ we have
\begin{equation}\label{eq:rightaction}
    r(X) f : g \mapsto
             \partial_{t}f(g\,e^{tX})|_{t=0}
             = (X_{L}f)(g),
\end{equation}
where $X_{L}$ is the left-invariant vector field taking value $X$
at the identity. Similarly, we obtain $l(X)f = X_{R}f$, where
$X_{R}$ is the corresponding right-invariant vector field.

Given an $\h$-representation $M$ with $\h$-action $m$, we view the
Hilbert space $L^{2}(G)\otimes M$ of sections of the trivial
bundle $G\times M$ as a representation of $\g\oplus\h$, taking the
left $\g$-action $l\otimes 1$ on the $L^{2}(G)$ component and the
diagonal right $\h$-action $r\otimes 1 + 1\otimes m$. This
apparent asymmetry between the left and right actions is a
consequence of the convention of trivializing vector bundles by
left translation. We say that a linear differential operator
$D:L^{2}(G)\otimes M\rightarrow L^{2}(G)\otimes N$ on $G$ is
\emph{homogeneous} relative to $\h$ if it commutes with the
$\g\oplus\h$ actions on its domain and range.

If the $\h$-actions on $M$ and $N$ exponentiate to give
single-valued actions of the Lie group $H$, then the
representations $M$ and $N$ induce equivariant $G$-bundles
$G\times_{H}M$ and $G\times_{H}N$ over the homogeneous space
$G/H$. The sections of these bundles correspond to the right
$H$-equivariant functions on $G$ taking values in $M$ and $N$
respectively, giving
\begin{equation*}
    L^{2}(G\times_{H}M) \cong (L^{2}(G)\otimes M)^{H}, \quad
    L^{2}(G\times_{H}N) \cong (L^{2}(G)\otimes N)^{H}.
\end{equation*}
If $D$ is homogeneous relative to $H$, then it preserves the
$H$-invariance of both its domain and range, and so it restricts
to give a $G$-equivariant operator
\begin{equation*}
    D_{0} : L^{2}(G\times_{H}M) \rightarrow L^{2}(G\times_{H}N).
\end{equation*}
In general, such a $G$-equivariant linear differential operator on
$G/H$ is called a \emph{homogeneous differential operator}.
However, the notion of homogeneous operators given above is more
flexible, as it is not limited to bundles that descend to the
quotent $G/H$, and we are free to work upstairs on $G$. In fact, a
homogeneous operator on $G$ defines not just one but rather an
entire family of twisted homogeneous operators on $G/H$ indexed by
the irreducible representations $U_{\mu}$ of $H$, obtained by
restricting to the subspaces of the domain and range transforming
like the dual representations $U_{\mu}^{\ast}$. Specifically, the
family of operators $D_{\mu}$ are given by
\begin{align*}
    D_{\mu} : \Hom_{\h}\bigl(U_{\mu}^{\ast},L^{2}(G)\otimes M\bigr)
              & \longrightarrow
              \Hom_{\h}\bigl(U_{\mu}^{\ast},L^{2}(G)\otimes N\bigr) \\
              L^{2}\bigl(G\times_{H}(M\otimes U_{\mu})\bigr)
              & \longrightarrow
              L^{2}\bigl(G\times_{H}(N\otimes
              U_{\mu})\bigr).
\end{align*}
Also, even if the $\h$-actions on $M$ and $N$ fail to exponentiate
to the Lie group $H$, the operators $D_{\mu}$ nevertheless descend
to well-defined operators on $G/H$ provided that the tensor
products $M\otimes U_{\mu}$ and $N\otimes U_{\mu}$ are indeed true
representations of $H$. In other words, twisting by a projective
representation with the opposite cocycle kills the obstruction.

In \cite{B}, Bott showed that the index of a homogeneous
differential operator, like that of a finite dimensional operator,
depends only on the domain and range and not on the operator
itself. Since the domain and range are now infinite dimensional
representations of $G$, Bott viewed them as elements of the
completion $\hat{R}(G)$ of the representation ring of $G$,
consisting of all possibly infinite formal sums
$\sum_{\lambda}a_{\lambda}[V_{\lambda}]$ with integer coefficients
indexed by the equivalence classes of irreducible representations
of $G$. In this notation, the space of sections of a homogeneous
vector bundle induces the class
\begin{equation}\label{eq:frobenius}\begin{split}
    \bigl[L^{2}(G\times_{H}M)\bigr] & = \sum_{\lambda}\,
    [V_{\lambda}]
    \dim\Hom_{G}\bigl(V_{\lambda},L^{2}(G\times_{H}M)\bigr) \\
    & = \sum_{\lambda}\,
    [V_{\lambda}]
    \dim\Hom_{H}( V_{\lambda},M ),
\end{split}\end{equation}
where the second line follows from the first by Frobenius
reciprocity.

\begin{theorem}[Bott]\label{th:bott}
    If $D:L^{2}(G\times_{H}M)\rightarrow L^{2}(G\times_{H}N)$ is
    an elliptic homogeneous differential operator on $G/H$, then
    its $G$-index is the element of $\hat{R}(G)$ given by
    \begin{equation*}
        \Index_{G}D = [L^{2}(G\times_{H}M)] - [L^{2}(G\times_{H}N)].
    \end{equation*}
    Furthermore, this difference is actually a finite element
    in $R(G)\subset \hat{R}(G)$.
\end{theorem}

Bott's theorem follows from the Peter-Weyl decomposition of the
space of $L^{2}$ functions on $G$ into the Hilbert space direct
sum
\begin{equation*}
    L^{2}(G) \cong \widehat{\bigoplus}_{\lambda} V_{\lambda} \otimes
    V_{\lambda}^{\ast}
\end{equation*}
with respect to the natural action (\ref{eq:ggaction}) of $G\times
G$ on functions. For sections of a homogeneous vector bundle, the
Peter-Weyl decomposition becomes
\begin{equation}\label{eq:peter-weyl}\begin{split}
    L^{2}(G\times_{H}M) & \cong \widehat{\bigoplus}_{\lambda}
    V_{\lambda} \otimes ( V_{\lambda}^{\ast} \otimes M )^{H} \\
    & \cong \widehat{\bigoplus}_{\lambda}
    V_{\lambda} \otimes \Hom_{H}( V_{\lambda},M ),
\end{split}\end{equation}
which shows that the expression (\ref{eq:frobenius}) for the class
of $L^{2}(G\times_{H}M)$ in $\hat{R}(G)$ completely characterizes
this space of sections as a representation of $G$. Since the
operator $D$ is equivariant with respect to the $G$-actions on its
domain and range, it can be written in block diagonal form as $D =
\bigoplus_{\lambda} D|_{V_{\lambda}}$, where each of the operators
\begin{equation*}
    D|_{V_{\lambda}}:V_{\lambda}\otimes\Hom_{H}(V_{\lambda},M)
    \rightarrow      V_{\lambda}\otimes\Hom_{H}(V_{\lambda},N)
\end{equation*}
is finite dimensional and must therefore have $G$-index
\begin{equation*}
    \Index_{G} D|_{V_{\lambda}} = [V_{\lambda}]\,\bigl(
    \dim\Hom_{H}(V_{\lambda},M) - \dim\Hom_{H}(V_{\lambda},N)
    \bigr).
\end{equation*}
The total $G$-index of $D$ is then the sum $\Index_{G}D =
\sum_{\lambda}\Index_{G}D|_{V_{\lambda}}$, and all but finitely
many of these summands must vanish since $D$ is Fredholm.

\section{The Geometric Dirac Operator}
\label{section:newstuff}

Returning to Kostant's Dirac operator (\ref{eq:kostantdirac}),
when applied to spinors $L^{2}(G)\otimes\bbS$ with $\g$ acting on
$L^{2}(G)$ by the right action (\ref{eq:rightaction}), it becomes
the operator
\begin{equation}\label{eq:geodiracoperator}
    \dirac = \sum_{i} c(X_{i}^{\ast})\,r(X_{i}) + c(v),
\end{equation}
where $c$ is the Clifford action on the spin representation
$\bbS$, and $v\in\Cl(\p)$ corresponds to the fundamental 3-form
(\ref{eq:threeform}). Recalling that the infinitesimal right
action of $\g$ on functions is the same as differentiation with
respect to the left-invariant vector fields, we see that this
Dirac operator is a linear differential operator. Furthermore,
since the expression (\ref{eq:geodiracoperator}) for the Dirac
operator is written entirely in terms of the right action on
$L^{2}(G)$ and assorted endomorphims of the $\bbS$ component, it
must automatically commute with the left action $l\otimes 1$ of
$\g$ on $L^{2}(G)\otimes\bbS$. A quick computation at the level of
the non-abelian Weil algebra then shows that $\dirac$ commutes
with the diagonal action $r\otimes 1 + 1\otimes \ads$ of $\h$ on
$V\otimes\bbS$ for any representation $V$ of $\g$. Therefore, the
Dirac operator $\dirac : L^{2}(G)\otimes\bbS\rightarrow
L^{2}(G)\otimes\bbS$ is homogeneous with respect to $\h$.

For reasons that will soon become evident, rather than decomposing
$\bbS$ into the two half-spin representations $\bbS^{+}$ and
$\bbS^{-}$ as usual, we instead take their dual representations
$\bbS_{+}^{\ast}$ and $\bbS_{-}^{\ast}$. This modification simply
introduces an overall sign factor, since the total spin
representation itself is self-dual. As for the half-spin
representations, they are self-dual when $\frac{1}{2}\dim\p$ is
even, and they are dual to each other when $\frac{1}{2}\dim\p$ is
odd.

\begin{theorem}\label{th:dirac-index}
Given an irreducible representation $U_{\mu}$ of $\h$ with highest
weight $\mu$ such that the tensor product $\bbS\otimes U_{\mu}$ is
a true representation of $H$, then the $G$-equivariant index of
the twisted Dirac operator
\begin{equation*}
    \dirac_{\mu}:
    L^{2}\bigl(G\times_{H}(\bbS_{+}^{\ast}\otimes U_{\mu})\bigr)
    \rightarrow
    L^{2}\bigl(G\times_{H}(\bbS_{-}^{\ast}\otimes U_{\mu})\bigr)
\end{equation*}
is $\Index_{G}\dirac_{\mu} =
(-1)^{w}[V_{w(\mu+\rho_{H})-\rho_{G}}]$ if there exists a Weyl
group element $w\in W_{G}$ such that the weight
$w(\mu+\rho_{H})-\rho_{G}$ is dominant for $G$. If no such element
$w$ exists, then $\Index_{G}\dirac_{\mu} = 0$.
\end{theorem}

\begin{proof}
    For such a choice of $U_{\mu}$, the operator $\dirac_{\mu}$
    descends to give a homogeneous differential operator on $G/H$.
    The symbol of the Dirac operator is Clifford multiplication,
    which is invertible, so the operator is elliptic. We may
    therefore apply Theorem \ref{th:bott} to compute its
    index, and by (\ref{eq:frobenius}) we have
    \begin{equation*}\begin{split}
        \Index_{G}\dirac_{\mu}
        & = \sum_{\lambda}\,[V_{\lambda}]
        \bigl(
            \dim\Hom_{H}(V_{\lambda},\,\bbS_{+}^{\ast}\otimes U_{\mu}) -
            \dim\Hom_{H}(V_{\lambda},\,\bbS_{-}^{\ast}\otimes U_{\mu})
        \bigr) \\
        & = \sum_{\lambda}\,[V_{\lambda}]
        \bigl(
            \dim\Hom_{H}(V_{\lambda}\otimes \bbS^{+},\,U_{\mu}) -
            \dim\Hom_{H}(V_{\lambda}\otimes \bbS^{-},\,U_{\mu})
        \bigr) \\
        & = \sum_{\lambda}\sum_{c\in C}
            (-1)^{c}\,[V_{\lambda}]\,\delta_{c\bullet\lambda,\,\mu},
    \end{split}\end{equation*}
    using the homogeneous Weyl formula (\ref{eq:homogeneous-weyl})
    in the final line to decompose the virtual representation
    $V_{\lambda}\otimes\bbS^{+}-V_{\lambda}\otimes\bbS^{-}$
    into irreducible representations of $\h$. Finally, we have
    $\mu = c\bullet\lambda = c(\lambda+\rho_{G})-\rho_{H}$ if
    and only if $\lambda = w(\mu + \rho_{H}) - \rho_{G}$ for
    $w = c^{-1}$.
\end{proof}

Note that even if $\bbS\otimes U_{\mu}$ is only a projective
representation of $H$, rather than a true representation, the
operator $\dirac_{\mu}$ is nevertheless Fredholm and
$G$-equivariant. Although it no longer descends to give an
elliptic operator on $G/H$, the proofs of both Theorem
\ref{th:bott} and Theorem \ref{th:dirac-index} continue to hold,
and the statement of Theorem \ref{th:dirac-index} for the
$G$-index of $\dirac_{\mu}$ is unchanged.

By taking the index of the Dirac operator on $G/H$ twisted by an
irreducible representation $U_{\mu}$ of $\h$, we have effectively
inverted the construction of Gross, Kostant, Ramond, and
Sternberg. Instead of using the homogeneous Weyl formula to
extract a multiplet of $\h$-representations associated to a given
representation of $\g$, we can now start with a single
representation of $\h$ and use this index to determine the unique
$\g$-representation from which it can be obtained. In fact, we can
be even more precise. In Theorem \ref{th:dirac-index}, the index
of the Dirac operator depends only on the domain and range and not
on the operator itself. However, we have been using a particular
choice of Dirac operator, which Kostant constructed in \cite{K}
specifically for the properties of its kernel and cokernel. Using
Kostant's results, we obtain a short proof of the following
theorem of Slebarski (see \cite{Sl}):

\begin{theorem}[Slebarski]\label{th:slebarski}
    Given an $\h$-representation $U_{\mu}$ as in the statement
    of Theorem \ref{th:dirac-index}, the space of harmonic spinors
    for the twisted Dirac operator
    \begin{equation*}
        \dirac_{\mu}:
        L^{2}\bigl(G\times_{H}(\bbS\otimes U_{\mu})\bigr)
        \rightarrow
        L^{2}\bigl(G\times_{H}(\bbS\otimes U_{\mu})\bigr)
\end{equation*}
    is $\Ker\dirac_{\mu} = V_{w(\mu+\rho_{H})-\rho_{G}}$ if
    there exists a Weyl element $w\in W_{G}$ satisfying
    the conditions of Theorem \ref{th:dirac-index}, and
    $\Ker\dirac_{\mu}=0$ otherwise.
\end{theorem}

\begin{proof}
    As in the proof of Theorem \ref{th:bott}, we use the homogeneous
    form (\ref{eq:peter-weyl}) of the Peter-Weyl theorem to
    decompose the kernel of $\dirac_{\mu}$ as the direct sum
    \begin{equation*}
        \Ker\dirac_{\mu} =
        \widehat{\bigoplus}_{\lambda}\Ker\dirac_{\mu}|_{V_{\lambda}},
    \end{equation*}
    over the finite dimensional Dirac operators
    $\dirac_{\mu}|_{V_{\lambda}}$ acting on the spaces
    \begin{equation*}
        V_{\lambda}\otimes\Hom_{H}(V_{\lambda},\,\bbS\otimes U_{\mu})
        \cong
        V_{\lambda}\otimes\Hom_{H}(U_{\mu},\,V_{\lambda}\otimes\bbS),
    \end{equation*}
    where all but finitely many of these kernels vanish since
    the operator $\dirac_{\mu}$ is Fredholm. Ignoring the signs of the
    half-spin representations, equation (\ref{eq:kostant-kernel})
    for the kernel of $\dirac$ implies that the kernel of
    $\dirac_{\mu}|_{V_{\lambda}}$ is $V_{\lambda}$ if $\mu =
    c\bullet\lambda$ for some $c\in C$ and $0$ otherwise.
\end{proof}

This result provides an explicit construction for any
representation $V_{\lambda}$ of $G$ as a space of twisted harmonic
spinors on $G/H$, giving a harmonic induction map from the
irreducible representations of $H$ to the irreducible
representations of $G$. In particular, if $H$ is a maximal torus
in $G$, then this theorem becomes a version of the Borel-Weil-Bott
theorem expressed not in its customary form involving holomorphic
sections and Dolbeault cohomology, but rather in terms of spinors
and the Dirac operator. Again, note that if $\bbS\otimes U_{\mu}$
is not a true representation of $H$, the statement and proof of
Theorem \ref{th:slebarski} still hold, but the twisted spinors and
Dirac operator no longer descend to the homogeneous space $G/H$.
Also note that we can recover the sign factor $(-1)^{c} =
(-1)^{w}$ in the homogeneous Weyl
formula~(\ref{eq:homogeneous-weyl}) and the index of the twisted
Dirac operator by splitting the spin representation into its dual
half-spin representations as we did in
Theorem~\ref{th:dirac-index}. The representations with positive
sign then appear in the kernel of the Dirac operator, while the
negative ones appear in the kernel of its adjoint, as in
(\ref{eq:kostant-kernel}).

\end{document}